\documentclass[10pt,a4paper]{amsart}
\usepackage{amssymb,latexsym,amscd,amsthm,amsfonts,verbatim,amsopn,color,amscd}

\usepackage{psfrag}
\usepackage[dvips]{epsfig}
\usepackage{color}

\usepackage[all,2cell,ps]{xy}
\usepackage[all]{xy}

\newtheorem{thm}{Theorem}[section]
\newtheorem{df} [thm]{Definition}
\newtheorem{lm}  [thm]{Lemma}
\newtheorem{crl} [thm]{Corollary}
\newtheorem{prop}[thm]{Proposition}

\newtheorem{rem}[thm]{Remark}

\newtheorem{expl}[thm]{Example}

\newcommand{\us}{\underline{\sigma}}
\newcommand{\Dl}{\mathcal{Q}_{n}(G)}
\newcommand{\Pno}{\Pi_{n, 0}}
\newcommand{\dl}{\mathcal{Q}_n^0(G)}
\newcommand{\uo}{\underline{\omega}}

\newcommand{\nG}{[n]\times G}

\newcommand{\TG}{\mathcal{T}_n^G}
\newcommand{\orb}{\mathfrak{o}}
\newcommand{\II}{\mathcal{I}}
\newcommand{\JJ}{\mathcal{J}}
\newcommand{\NN}{\mathcal{N}}
\newcommand{\LL}{\mathcal{L}}
\newcommand{\GG}{\mathcal{G}}

\newcommand{\KK}{\mathcal{K}}

\newcommand{\TT}{\mathcal{T}}
\newcommand{\TNG}{\mathcal{T}_n(G)}
\newcommand{\RD}{\widetilde{\Delta}}
\newcommand{\RN}{\widetilde{\NN}}

\newenvironment{pf}{\noindent {\bf Proof. }}{\hfill $\Box$\vspace{0.3cm}}

\title[Dowling trees]{Nested set complexes of Dowling lattices\\ and complexes of Dowling trees}

\author{Emanuele Delucchi}
\address{Dipartimento di matematica, Università di Pisa, largo Bruno
  Pontecorvo 5, Pisa, Italy}
\email{delucchi@mail.dm.unipi.it}
\thanks {Research partially supported by the Swiss National Science
  Foundation, project PP002--106403/1}

\begin{document}

\maketitle

\begin{abstract}Given a finite group $G$ and a natural number $n$, we study
  the structure of the complex of nested sets of the associated
  Dowling lattice $\Dl$ (see \cite{Dow}) and of  its subposet of the
  $G$-symmetric partitions $\dl$ which was recently introduced by
  Hultman in \cite{Hu}, together with the complex of $G$-symmetric
  phylogenetic trees $\TG$. Hultman shows that the complexes $\TG$ and $\RD(\dl)$ are homotopy
  equivalent and Cohen-Macaulay, and determines the rank of their top
  homology. 

An application of the theory of building sets and nested set complexes by Feichtner and
 Kozlov \cite{FK} shows that in fact $\TG$ is subdivided by the order complex
 of  $\dl$. We introduce the complex of Dowling trees $\TNG$ and prove
 that it is subdivided by the order complex of $\Dl$ and contains
 $\TG$ as a subcomplex. We show that $\TNG$ is obtained from $\TG$ by
 successive coning over certain subcomplexes.
It is well known that $\Dl$ is shellable, and of the same dimension as $\TG$.
 We explicitly and independently calculate how many homology spheres are added in passing from $\TG$ to  $\TNG$.\end{abstract}

\section*{Introduction}

Dowling lattices are named after T.\ A. Dowling, who first studied a particular class of arrangements of hyperplanes 
whose intersection lattices can be obtained by enriching the partition
lattice with sets of elements of the cyclic groups $\mathbb{Z}_n$ \cite{Dow1}. In a separate work (\cite{Dow}), Dowling took the step to a general combinatorial point of view, and thoroughly studied what he called `posets based on finite groups'. We wish to point to \cite{Dow} as a very readable and comprehensive introduction to the subject. \\
Both approaches to these structures were followed since the work of
Dowling. Let us mention, as two examples, the work of Ehrenborg and
Readdy \cite{ER, ER2}, who introduced the notion of {\em Dowling transform}
of an arrangement of hyperplanes and showed that this transformation
preserves supersolvability, and the study of a combinatorial
generalization of Dowling lattices by Hanlon \cite{Han}.\\

The theory of building sets and nested set complexes was initiated and
developed by Feichtner and Kozlov in \cite{FK} as the combinatorial
framework of the De Concini - Procesi models for hyperplane
arrangements. A complex of {\em nested sets}
is associated to any meet-semilattice and any of its {\em building
  sets} - i.e., subsets of the semilattice satisfying some
conditions that are inspired by the special properties of the set
of irreducible elements in a geometric lattice.\\
The study of such structures was carried out in different
contexts (see \cite{FK2, FY, FM, F, D, FS, CD}), leading to
new results or sharpening the understanding of previously studied
objects. In particular, this theory has proved to be a very useful
tool for the study of {\em complexes of trees} of various kind (see
\cite{F, D}).\\

The study of abstract simplicial complexes whose cells are indexed by
combinatorial types of rooted trees on a fixed number of leaves was
recently brought to a broad attention by the work of Billera, Holmes and Vogtmann \cite{BHV}. They considered
the space of all possible phylogenetic trees of a certain set of
biological species, with lengths on edges representing the genetic
distance of two mutations. This space is a cone with apex the unique
tree with all edge lengths equal to zero. The base of this cone (the
link of the apex) is obtained by considering the trees with unit edge
lengths (i.e., cutting the cone by a transversal hyperplane). This
space has a natural stratification, with each cell corresponding to a
combinatorial type of trees - this is the {\em complex of phylogenetic
trees.}\\
Complexes of trees were studied for different purposes, before and
after \cite{BHV}, see \cite{Bo, Han1, HW, RW, V, TZ}.  In particular,
in \cite{RW} Robinson
and Whitehouse determined the homotopy type of the complex of
phylogenetic trees on $n$ leaves to be a wedge of $(n-1)!$ spheres
of dimension $n-3$. Later on, shellability of these complexes was proved by Trappmann
and Ziegler \cite{TZ} and, independently, by Wachs
(unpublished, acknowledged in \cite{F}). 
Ardila and Klivans \cite{AK} proved that the complex of trees can be subdivided
by the order complex of the partition lattice.
 One of the recent
applications of the theory of nested set complexes sharpened this last
result: in \cite{F}, Feichtner showed
that the order complex of the partition lattice is obtained from the
complex of trees by a sequence of stellar subdivisions. This property was
proven to hold even in the more general class of $k$-trees (see
\cite{D}).\\

Our purpose is to study the nested set complexes of Dowling lattices
and look for a connection with the complex of trees.
The initial motivation came from the work of Hultman \cite{Hu}, who defined the
{\em complex of $G$-symmetric phylogenetic trees} and showed that it
has the same homotopy type as a particular subposet of the Dowling
lattice. We sharpen this homotopy equivalence by showing that in fact the two complexes are related by a sequence of
stellar subdivisions (Theorem \ref{subd}) and are therefore homeomorphic. Moreover, we put these
objects into the classical theory of Dowling lattices: we introduce a
combinatorially and topologically suitable notion of Dowling trees
(Definition \ref{Dtrees}, Corollary \ref{subD}) and study their
topological relationship with Hultman's complexes (Theorem \ref{nostro}, Remark \ref{ultimo}).\\

This paper is organized as follows: in sections \ref{pbofg} and
\ref{section_trees} we give a detailed picture
of Dowling lattices and $G$-symmetric phylogenetic trees, reviewing
the basics and developing a notation that will prove to be appropriate
for a direct application of the theory of nested set complexes.
They will enter the picture in Section \ref{subds}, where the
result of Hultman is sharpened by showing that the complex of
$G$-symmetric phylogenetic trees actually is subdivided by the associated
subposet of the Dowling lattice. 
In Section \ref{section_dtrees} the notion of a {\em Dowling
  tree} is introduced as naturally associated to the nested set
complex of the full Dowling lattice, thus being probably the
appropriate Dowling generalization of the complex of phylogenetic
trees.
The complex of Dowling trees contains the complex of $G$-symmetric
phylogenetic trees as a simplicial subcomplex. In the last section we
describe explicitly  how the bigger complex can be obtained from the
smaller one by successively coning over certain subcomplexes. By
keeping under control the topology of those subcomplexes, we calculate
how many homology spheres arise at each step, thus explicitly
relating the ``numerologically" suggestive expressions of the top
homology ranks of the two complexes.

\section{Posets based on finite groups}\label{pbofg}

Before starting out, let us fix some general notations. In this paper we will deal with finite partially ordered sets, briefly called {\em posets}. The main topological structure associated to a poset $P$ is its {\em order complex} $\Delta(P)$, the (abstract) simplicial complex of the totally ordered subsets of $P$. We will almost never distinguish between an abstract simplicial complex and its geometric realization, thus speaking of `topological properties' of an abstract simplicial complex. It is easily seen that if $P$ possesses a maximal element (that is customarily denoted $\hat{1}$), then $\Delta(P)$ is a cone over $\Delta(P\setminus\{\hat{1}\})$. An analogous statement holds of course if $P$ has a minimal element (usually denoted $\hat{0}$). To capture the `essential' topological information we define the {\em reduced order complex} $\RD(P)$ as the order complex of the poset obtained from $P$ by removing the maximal and the minimal elements, if $P$ has any.

In considering partitions, we will switch between the set-theoretic
notation $\sigma=S_0\amalg S_1\amalg \dots \amalg S_k$ and the (more
customary) `block notation' $\sigma=S_0\vert\dots\vert S_k$. Sets of
partitions can be ordered {\em by refinement}, i.e., setting
$\sigma'<\sigma $ if every block of $\sigma'$ is contained in a block
of $\sigma$.

\subsection{Dowling lattices} 

{\df Let $G$ be a finite group and $n$ a natural number. Consider the
  action of $G$ on the set $\{0\}\cup ([n] \times G)$ defined by
  $g((i,h))=(i,gh)$ and $g(0)=0$. The partitions of $\{0\}\cup ([n]
  \times G)$ such that this action induces an action on the blocks are
  called {\bf\em $\mathbf{G}$-symmetric partitions}. A block of a
  $G$-symmetric partition is called {\bf\em simple} if its orbit under
  this action has length $|G|$. 

We call $\Dl$ the set of $G$-symmetric partitions such that the only non-simple block is the block containing $0$. The ordering by refinement turns it into a lattice, called the {\bf\em Dowling lattice}.}\\

The `forgetful' map $\{0\}\cup ([n] \times G) \rightarrow \{0\} \cup [n]$ defined by $(i,g)\mapsto i$ and $0\mapsto 0$ induces a mapping of $\Dl$ onto the poset $\Pno$ of partitions of the set $\{0\}\cup [n]$. This map sends $\omega\in \Dl$ to its {\em associated partition} $\uo \in \Pno\cong \Pi_{n+1}$. An element of $\Pno$ will be written as $\alpha:= A_0\vert A_1 \vert \dots \vert A_k$, where we agree to choose the indexing such that $0$ is always contained in the block indexed by $0$.
Of course, $\alpha=\uo$ has one block $A_i$ for every orbit of blocks in $\omega$, and the $A_i$ with $i>1$ correspond to the simple blocks.

We now see that we can encode in a unique way any $\omega \in \Dl$ in the following data:
\begin{itemize}
\item[$\bullet$] A partition $\alpha\in\Pno$, called the {\em associated
  partition} of $\omega$.
\item[$\bullet$] For every block $A_j=\{i_1,\dots , i_k\}$ of $\alpha$ with $j>0$ an $|A_j|$-tuple\\
  $A_j^G:=(id,g^{\alpha}_{i_2},g^{\alpha}_{i_3}, \dots , g^{\alpha}_{i_k})\subset G^{k}$.
\end{itemize} 

This gives a good encoding of $G$-symmetric partitions if we agree to always order the elements in the blocks $A_j$ in increasing order: $i_1< i_2 <\dots < i_k$.

The order relation translates to the following: $\omega_1 <
\omega_2$ if and only if
\begin{itemize}\item[$\bullet$] the associated partitions satisfy $\uo_1<\uo_2$ \item[$\bullet$] for every block $B$ of $\uo_2$, if $C_1,\dots ,C_k$ are the
blocks of $\uo_1$ that subdivide it, there are elements $h_1,\dots
,h_k$ in $G$ such that for every $j\in B\cap C_i$ we have $g^{\uo_2}_j=h_i
g^{\uo_1}_j$. \end{itemize}

Dowling gave a complete and very readable survey on these objects in
his seminal paper \cite{Dow}, to which we point as a general
reference. Although the definition we gave is inspired by the setting of
\cite{Hu}, it can be easily seen to be equivalent to the definition
of Dowling, e.g. by thinking in terms of the encoding we presented above.

The main statement on the topology of these lattices that can be deduced from \cite{Dow} is summarized in the following proposition.

{\prop The reduced order complex $\RD(\Dl)$ is homotopy equivalent to a wedge of $(|G|+1)(2|G|+1)\dots ((n-1)|G|+1)$ spheres of dimension $(n-2)$. }\\

{\pf The proof is a concatenation of arguments from \cite{Dow} (where
    $\Dl$ is shown to be supersolvable),
    \cite{Bj} (where it is proved that supersolvable lattices are
    shellable), and \cite{Sta} (for the enumeration of the number of
    homology spheres). Details are left to the reader. \hfill$\square$}

\subsection{The subposet $\mathbf{\dl}$} In \cite{Hu}, the author
introduces a subposet of $\Dl$ that is of particular interest in
connection with $G$-symmetric phylogenetic trees.

{\df We define $\dl$ to be the subposet of $\Dl$ consisting of the
  partitions with trivial zero block.}

{\rem Our definition is slightly different from the one given in
  \cite{Hu}, because we allow the minimal element of $\Dl$ to be in $\dl$.}\\

In particular, given $\sigma \in \dl$, every element of $\nG$ is
contained in a simple block of $\sigma$.

In general, $\dl$ is not a lattice. Indeed, consider the following two
elements of $\mathcal{Q}_3^0(\mathbb{Z}_2)$:
\begin{center}
$\sigma_1:= \{0\}\{(1,0),(2,0)\}\{(1,1)(2,1)\}\{(3,0)\}\{(3,1)\}$,\\
$\sigma_2:= \{0\}\{(1,0),(2,1)\}\{(1,1)(2,0)\}\{(3,0)\}\{(3,1)\}$.
\end{center}
Their join is not contained in $\mathcal{Q}_3^{\mathbb{Z}_2}$, where
there is no element that is bigger of both $\sigma_1$ and $\sigma_2$.\\

Nevertheless, $\dl$ is a meet-semilattice. We prove the
following easy lemma.

{\lm \label{lower} For any $\sigma \in \dl$, we have an isomorphism $$(\dl)_{\leq \,\sigma} \simeq (\Pi_n)_{\leq \,\us }.$$ }
{\pf The map defined by $\sigma \mapsto \us$ is clearly a poset morphism,
and obviously surjective. For injectivity note that, given $\us'\leq
\us$, any set of representatives of the orbit blocks for $\sigma$ forces
the choice of the $\ell$-tuples associated to the blocks to
$\us'$. Indeed, if $B$ is a block of $\us$ that is associated to $(id,
g_2^{\sigma}, \dots , g_\ell^{\sigma})$, then if $C=\{i_1, \dots ,
i_k\}\subseteq B$ is a block of $\us'$, the only possibility for a
partition associated to $\us'$ to be below $\sigma$ is to associate to
$C$ the $k$-tuple $(id, g_{i_2}^{\sigma}(g_{i_1}^{\sigma})^{-1}, \dots , g_{i_k}^{\sigma}(g^{\sigma}_{i_1})^{-1})$. \hfill $\square$ }

This makes the work of Section \ref{subds} possible, where the theory
of building sets and nested set complexes will be applied to the
posets $\dl$.\\

\newcommand{\edge}[1]{\ar@{-}[#1]}
\newcommand{\edgeh}[1]{\ar@{--}[#1]}
\newcommand{\node}{*+[o][F-]{ }}
\newcommand{\w}{\overline{1}}
\newcommand{\ww}{\overline{2}}
\newcommand{\www}{\overline{3}}

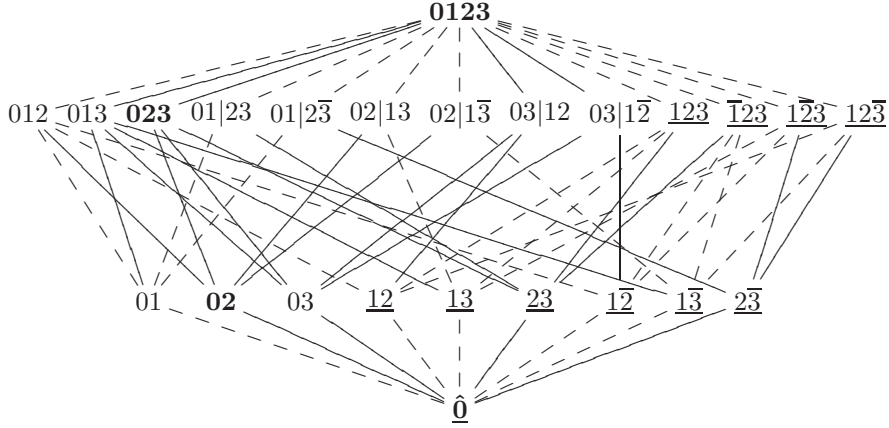
\begin{figure}[h]\label{Q3}
\centerline{
{\xymatrix@C-23pt{
&&&&&&\mathbf{0123}\edgeh{dllllll}\edge{dlllll}\edge{dllll}\edgeh{dlll}\edgeh{dll}\edgeh{dl}\edgeh{d}\edge{dr}\edge{drr}\edgeh{drrr}\edgeh{drrrr}\edgeh{drrrrr}\edgeh{drrrrrr}&&&& \\
012\edgeh{ddrr}\edge{ddrrr}\edgeh{ddrrrrr}\edgeh{ddrrrrrrrr}
                               &013\edge{ddr}\edge{ddrrr}\edge{ddrrrrr}\edge{ddrrrrrrrr}
                                                                &\mathbf{023}\edge{ddrr}\edge{ddr}\edge{ddrrrrr}
                                                                                                &01|23\edgeh{ddl}\edge{ddrrrr}
                                                                                                &01|2\www\edgeh{ddll}\edge{ddrrrrrr}
                                                                                                &02|13\edge{ddll}\edgeh{ddr}
                                                                                                &02|1\www\edge{ddlll}\edgeh{ddrrr}
                                                                                                &03|12\edge{ddlll}\edge{ddll}
                                                                                                &03|1\ww\edge{ddllll}\edge{dd}
                                                                                                &\underline{123}\edgeh{ddllll}\edgeh{ddlll}\edge{ddll}
                                                                                                                        &\underline{\w 23}\edge{ddlll}\edgeh{ddll}\edgeh{ddl}
                                                                                                                             &\underline{1\ww 3}\edgeh{ddlllll}\edgeh{ddlll}\edge{ddl} 
                                                                                                                                  & \underline{12\www}\edgeh{ddlllllll}\edgeh{ddlll}\edge{ddll}&  \\
&&&&&&&&&&&&&&\\
&& 01                            &\mathbf{02}                    &03                             &\underline{12}       &\underline{13}  
                                                                                                                             &\underline{23}  
                                                                                                                                  &\underline{1\ww}&\underline{1\www}&\underline{2\www}  \\
&&&&&&\underline{\mathbf{\hat{0}}}\edgeh{ullll}\edge{ulll}\edge{ull}\edgeh{ul}\edgeh{u}\edge{ur}\edgeh{urr}\edgeh{urrr}\edge{urrrr}&&&&
}}
}
\caption{The lattice $\mathcal{Q}_3(\mathbb{Z}_2)$ and its subposet $\mathcal{Q}_3^0(\mathbb{Z}_2)$ (underlined).
The bold elements give a modular chain, the dashed chains are the $15$ homology chains of the corresponding shelling.}
\end{figure}

As an example, we depict in Figure 1 
the Dowling lattice $\mathcal{Q}_3(\mathbb{Z}_2)$, where we write only one representative for the nonsingleton blocks. The numbers refer to the associated partitions, and an overline over an element indicates that this number is associated with the nonidentity element of $\mathbb{Z}_2$. 
Thus we will write $01|2\overline{3}$ instead of $\{0,(1,0), (1,1)\}\{(2,0),(3,1)\}\{(2,1),(3,0)\}$. 
The elements of $\mathcal{Q}_3^0(\mathbb{Z}_2)$ are underlined.

The homotopy type of $\RD(\dl)$ was determined in \cite{Hu} by comparison with the complex of $G$-symmetric phylogenetic trees (see Section \ref{section_trees}), as an application of discrete Morse theory.

\section{G-Symmetric phylogenetic trees}\label{section_trees}

Recall that we fixed once and for all a natural number $n$. In this
context, given a finite group $G$, a $G$-tree is a rooted
tree whose leaves are in bijection with the set $\nG$. The group $G$
acts on the set of leaves by means of the 'standard' action described
in the previous section.

We now specify a class of $G$-trees that will be the object of our
study.

{\df\label{hultman-tree} A {\bf \em $\mathbf{G}$-symmetric phylogenetic tree} is a
  $G$-tree satisfying the following additional conditions:
\begin{itemize}
\item[(1)] Every internal vertex (except the root) has degree at least $3$.
\item[(2)] The tree is invariant under the $G$-action.
\item[(3)] For any two different elements $g,h \in G$ and any $i\in
  [n]$, the (unique) shortest path connecting the leafs labelled
  $(i,g)$ and $(i,h)$ passes through the root.
\end{itemize}
The set of $G$-symmetric phylogenetic trees is denoted by $\TG$.}\\

This is definition 3.1 of \cite{Hu}, where some properties of those
trees are listed. Here we need only recall that every internal edge $t$ of a
tree $T\in\TG$ generates an orbit $\orb (t)$ (called {\em inner orbit}) of cardinality $|G|$ under the
action of $G$. Following \cite{Hu}, one can associate to $\orb (t)$
the partition $\pi(\orb (t))$ of $\nG$ obtained by putting in the same block all
labels of leaves that are in the same connected component of $T$ after
removing all the edges in $\orb(t)$, and adding $0$ to the block
corresponding to the component containing the root. 
We want to slightly modify this definition. 

\subsection{Some notation} Note that in all partitions of $\{0\}\cup [n]$ that are associated to some $\sigma \in \dl$ the element $0$ is alone in its block. Therefore we may reconstruct every such partition from the corresponding partition of $[n]$ by just adding the block $\{0\}$. We then agree that, for $\sigma\in \dl$, in this section we will consider $\us \in \Pi_n$.\\

Now consider a tree $T$ satisfying definition \ref{hultman-tree}. For any vertex
$v$ of $T$ let $\lambda_v$ denote the set of leaf labels such that the path
connecting them to the root traverses $v$. 
Let 
$$\lambda_{Gv}:= \bigcup_{g\in G}\lambda_{gv}$$
\noindent denote the set of labels of leaves that are separated from the root by a vertex of the form $gv$ for some $g\in G$.

Then, given an inner edge $t$ of $T$, let $\lambda_t:= \lambda_v$ where $v$ is the vertex of $t$ that is further from the root. For every internal edge $t$ we can then define a partition  $\sigma(t)\in \dl$ as
$$ \coprod_{(i,g) \notin\lambda_{Gt}}\{(i,g)\}\coprod_{g\in
  G} \lambda_{gt}.  $$
It is clear that $\lambda_{gt}=\{(i,gh)\vert (i,h)\in \lambda_t\}=g\lambda_t$, and
therefore we see that $\sigma(t)$ has only one
nonsingleton orbit. Thus the associated partition $\us(t)\in \Pi_n$
has only one nonsingleton block.


The next two subsections present some material of \cite{Hu} in a
language and from a viewpoint that are well-suited to our methods.

\subsection{Inner orbit contraction} The contraction of all edges in
an inner orbit $\orb (t)$ turns a $G$-symmetric tree $T$ into another tree
$T'\in \TG$. We can then define the following partial order on $\TG$.

{\df\label{order} Given $T,T'\in \TG$, define $T\leq T'$ if and only if $T'$ can be
obtained from $T$ by a sequence of inner orbit contractions.}\\

The importance of this ordering is shown in the following proposition.

{\prop {\bf \cite[Corollary 3.5]{Hu}} With the partial ordering of
  definition \ref{order}, $\TG$ is the face poset of a
  pure simplicial complex of dimension $n-2$. }\\

Before turning our attention to the operation of inner orbit {\em extension}, which is inverse to the contraction defined above, let us see what kind of relation one can draw between partitions associated to different contractions on the same tree.

So suppose again a tree $T\in\TG$ be given, and consider two inner
edges $t, t'$ of $T$ that are not in the same $G$-orbit. We have seen
that the associated partition $\us(t)$ (respectively $\us(t')$) has
exactly one nonsingleton block, say $B$ (resp. $B'$), associated to the unique nonsingleton orbit of $\sigma(t)$ ($\sigma(t')$), of which we consider a representative block $S$ ($S'$). 
Suppose we first contract the edges in the orbit $\orb (t)$.

If, for some $h\in G$, $S'\subseteq hS$, then by definition $gS'\subseteq g(hS)$ for all $g\in G$. We conclude that in this case $\sigma(t)>\sigma(t')$. Therefore, the unique nonsingleton block $B'$ of the partition $\us(t')$ is contained in $B$, so $\us(t')<\us(t)$.

If for some $h\in G$ the reverse inclusion $S'\supset hS$ holds, then of course the conclusion above holds with $t$ and $t'$ switched.

The fact that $T$ is a tree excludes the possibility that, if neither of the previous cases enters, $S'\cap hS \neq \emptyset$ for some $h\in G$.

We summarize the conclusion for later reference.

{\rem If $t, t'$ are two inner edges of a $G$-symmetric tree $T$ such
  that $\orb (t)\neq \orb (t')$ and if $\sigma(t)$ and $\sigma(t')$
  are incomparable, then the sets $\lambda_{Gt}$ and $\lambda_{Gt'}$
  are disjoint. In particular the associated nonsingleton blocks $B$,
  $B'$ of $\us(t)$ and $\us(t')$ are either contained in one
  another, or are disjoint.

\label{tree_nested} In particular, we can associate to every tree $T\in \TG$ a subset $N(T)\subset \dl$ defined as 
$$N(T):=\{\sigma(t) \mid t\textrm{ is (a representant of the orbit of) an inner edge of }T \}.$$ 
$N(T)$ has the property that the unique nonsingleton block orbits of any incomparable $\sigma, \sigma' \in N(T)$ are disjoint. }

\subsection{Inner orbit extension}

We now discuss the inverse of the above operation: {\em inner orbit extension}. 

For this, we suppose a tree $T\in \TG$ to be given together with a
partition $\sigma\in \dl$ that has exactly one nonsingleton orbit
$\orb$ (of which we consider a representative block $S$). The
preceding observations suggest to require the following condition to be satisfied by $\sigma$:
\begin{center}{$(\mathbf{\ast})\;$ }{\em For any inner edge $t$ of $T$, if neither $S \subset \lambda_t$ nor $S\supset \lambda_t$ then $S \cap \lambda_t = \emptyset$. }\end{center} 

In the following we will show how these data give rise to a tree $T'\in \TG$ such that $T$ is obtained from $T'$ by an inner orbit contraction that is represented by $\sigma$.

First of all it is clear that there is a unique vertex $v$ of $T$ such
that the component that is separated from the root by removing $v$ is
minimal with the property of containing all leaves labelled by elements
of $S$. The family of sets $\{ \lambda_{gv} \setminus
gS \vert g\in G\}$ can be obtained from the representative
$\lambda_v \setminus S$ by the action of $G$. 

Because of property ($\ast$), we may partition the edges $f$ incident to $v$ and with $\lambda_f\subset
\lambda_v$ into the two following classes:\begin{center}
$F:=\{f\vert\lambda_f \subset S\}$, $E:=\{f\vert\lambda_f\cap S =\emptyset  \}$.\end{center}

All is now prepared for the extension. We first delete all edges $f\in E\cup F$, and for each of them we get a connected component $T_f$ not containing the root.

We then grow an edge $t$ below $v$ appending all edges in $E$ (and
their whole components) to $v$, and the edges of $F$ (with their connected components) to the
other vertex of $t$.

In the tree that we have now constructed we clearly have $\lambda_t=S$. Of course we may repeat the whole process inserting edges
$gt$ below $gv$ for any $g\in G$, eventually reaching a
tree $T'$ satisfying the requirements.\\

The tree $T$ could have been reached by a sequence of inner orbit extensions. If we consider the set $N$ of all partitions that correspond to an edge orbit of $T$ we indeed may reconstruct $T$ starting from the unique tree without inner edges by recursively performing the above process with all elements of $N$.
We have seen that the unique condition enabling to perform such an extension is given at each step by ($\ast$). We summarize our considerations with the following two statements.\\

\noindent{\bf Condition N:} {\em Given a subset $X\subset \dl$ of partitions that have only one orbit consisting of non-singleton blocks, we say that $X$ {\em satisfies condition N} if for any two incomparable $\sigma, \sigma' \in X$ the only non-singleton blocks of $\us ,\us '$ are disjoint.}

{\rem To any set $N$ of one-nonsingleton-orbit partitions from $\dl$
  satisfying condition N we can naturally associate a tree $T$, and
  this is such that, with the notation of remark \ref{tree_nested}, $N=N(T)$. }

\section{Homeomorphism through subdivisions}\label{subds}

The reader familiar with the subject will have already noticed that
in the previous section all material has been prepared for a
direct application of the theory of building sets and nested set
complexes. This theory was first developed by Feichtner and Kozlov in
\cite{FK} as the combinatorial framework of the De Concini-Procesi
models for hyperplane arrangements. We refer to that paper for a
thorough introduction to this subject. Here we recall only the
main definitions.

{\df 

Let $\LL$ be a meet-semilattice. A {\bf building set} of $\LL$ is a subset $\GG\subseteq\LL\setminus\{\hat{0}\}$ such that for any $x\in\LL\setminus\hat{0}$ there is an isomorphism
\begin{displaymath}
\varphi_x: {\Large\prod_{\scriptstyle{{z\in\textrm{\textup{max}}}\,\GG_{\leq x}}}} [\hat{0},z] \rightarrow [\hat{0},x]
\end{displaymath}
with $\varphi_x(0,\dots, 0,z,0,\dots,0)=z$ for $z\in\textrm{\textup{max}}\,\GG_{\leq x}$.

We call a set $N\subseteq\GG$ {\bf\em nested} ($\GG$-nested, if specification is needed) if, for any set $\{x_1,\dots,x_\ell\}\subseteq N$ ($\ell\geq 2$) of incomparable elements, the join $x_1\vee \dots\vee x_\ell$ exists and is not an element of $\GG$.
The {\bf\em nested set complex} of $\LL$ with respect to $\GG$, denoted $\NN(\LL,\GG)$, is the abstract simplicial complex of all nonempty $\GG$-nested sets. If $\LL$ has a maximal element $\hat{1}$ and $\GG$ contains it, then the nested set complex is a cone with apex $\{\hat{1}\}$. The base of this cone is the {\bf\em reduced nested set complex} $\RN(\LL,\GG)$.}\\

One of the main topological features of this theory is the following
theorem which first appeared in \cite{FM} in a version for atomic lattices. It was then extended to its full generality in \cite{CD,D}, and to these papers we refer for a careful topological treatment of the concept of {\em stellar subdivision} of an abstract simplicial complex. Here we only mention that the geometric realizations of two abstract simplicial complexes that are related by subdivisions are homeomorphic (see \cite[Definition 2.1]{D}).

{\thm\label{subdiv} Consider two building sets $\GG_1,\GG_2$ in a meet-semilattice $\LL$. If $\GG_1\subset \GG_2$, then the simplicial complex $\NN(\LL, \GG_2)$ can be obtained from $\NN(\LL,\GG_1)$ by a sequence of stellar subdivisions.}\\

Note that, for any semilattice, there is a unique minimal building
set. It is given by the set of all elements $x$ such that
the interval $[\hat{0},x]$ cannot be decomposed in a product of
smaller principal order ideals. For example, in the partition lattice $\Pi_n$
those elements are the partitions with only one nonsingleton
block. The minimal building set of $\Pi_n$ will be denoted by $\II$.

On the other hand, for any meet-semilattice the maximal building set is the whole poset, and the associated reduced nested set complex is then $\RN(\LL, \LL)=\RD(\LL)$. We then have following corollary of the previous theorem.

{\crl \label{subd_crl} Let $\LL$ be a meet-semilattice and $\GG$ a building set in $\LL$. Then $\RD(\LL)$ can be obtained from $\RN(\LL,\GG)$ by a sequence of stellar subdivisions. }\\

In analogy with the partition poset let us define a subset $\II^G\subset \dl$ as
follows:
$$\II^G := \{\sigma\in \dl \vert \underline{\sigma}\in\II  \}. $$

The following proposition shows that this is indeed "the right
definition".

{\prop $\II^G$ is the minimal building set of $\dl$.}

{\pf With lemma \ref{lower} the claim follows immediately by comparison with
  $\Pi_n$.\hfill$\square$}\\

{\prop\label{tg} The complexes $\TG$ and $\NN (\II^G, \dl)$ are isomorphic.}

\nopagebreak

{\pf The condition of being nested in $\II^G$ is equivalent to
  condition N of the previous section. \hfill$\square$ }\\

We are ready to state  the main result of this section, which is now an easy
application of corollary \ref{subd_crl}.

{\thm \label{subd} The order complex $\RD(\dl)$ is obtained from the complex of
$G$-symmetric trees $\TG$ by a sequence of stellar subdivisions.}\\

Hultman calculated the homotopy type of $\TG$ in \cite{Hu}. We include this result in the following corollary, that is intended to summarize our topological knowledge about G-symmetric partitions and G-symmetric phylogenetic trees.

{\crl\label{hultman} The simplicial complexes $\RD(\dl)$ and $\TG$ are PL-homeomorphic. They are homotopy equivalent to a wedge of $$(|G|-1)(2|G|-1)\dots ((n-1)|G|-1)$$ spheres of dimension $(n-2)$.}

\section{Dowling Trees}\label{section_dtrees}

The natural task at this point is 
to study the nested set complexes of the full Dowling lattice $\Dl$.

\subsection{Nested set complexes in $\Dl$} We want to determine the minimal building set $\JJ^G$ of $\Dl$. By an easy check (or by comparing Theorem 2 of \cite{Dow}) one sees that, given any $\omega \in \Dl$ with zero block $S_0$ and fixed chosen orbit representatives $S_i$, $i=1,\dots , k$, there is a natural isomorphism 
$$(\Dl)_{\leq \omega} \stackrel{\sim}{\longrightarrow} \Pi^G_{m(\omega)}\times \Pi_{|S_1|}\times \dots \times \Pi_{|S_k|}$$
where $m(\omega):=\frac{(|S_0|-1)}{|G|}$.

With this decomposition, we see that any irreducible element having a nonzero block that is not a singleton must have $S_0=\{0\}$, thus be an element of $\II^G$.
Moreover, $\JJ^G \setminus \II^G$ consists of the partitions where all simple blocks are singletons. 

We will distinguish these two types of elements in $\JJ^G$ by calling any $x\in \JJ^G \cap \II^G$ {\em of type $1$}, while we will refer to the elements in $\JJ^G\setminus \II^G$ as to those {\em of type} $0$.\\

{\rem A subset $X$ of $\JJ^G$ is nested if and only if, for any  $\omega, \omega' \in X$, the only nonsingleton blocks of the associated partitions $\uo, \uo'$ of $\{0\}\cup [n]$ are either disjoint or contained in one another.}\\

The following facts are now at hand, and we collect them for later reference.

{\lm\label{minbuild} Let $\JJ^G$ denote the minimal building set of $\Dl$.\begin{itemize}
\item[(1)] $\JJ^G \cap \dl = \II^G$.
\item[(2)] $\NN(\II^G,\dl)\subseteq \NN(\JJ^G,\Dl)$
\item[(3)] For any $X\in \NN(\JJ^G,\Dl)$, $X\cap \II^G \in \NN(\II^G,\dl)$.
\item[(4)] If $X\in \NN(\JJ^G, \Dl)$, the subset $X\setminus \II^G$ of the elements of type $0$ is linearly ordered. 
\end{itemize}}

\subsection{Dowling trees}

The last question we want to address is whether the nested set complex of the full Dowling lattice has an interpretation in terms of trees. The answer is positive, and leads to the definition of what we would like to call {\bf \em the complex of Dowling trees}.

{\df\label{Dtrees} Given a natural number $n$ and a finite group $G$, a Dowling tree
  is a $G$-tree T with some distinguished vertices, called {\bf \em
    zero vertices}, satisfying the following conditions:
\begin{itemize}
\item[(0)] The root is a zero vertex.
\item[(1)] Every internal vertex (except the root) has degree at least $3$.
\item[(2)] The tree is invariant under the $G$-action, and the zero vertices are fixed by this action.
\item[(3)] For any two different elements $g,h \in G$ and any $i\in
  [n]$, the (unique) shortest path connecting the leafs labelled
  $(i,g)$ and $(i,h)$ passes through exactly one zero vertex.
\item[(4)] The zero vertices are the vertices of a path beginning at the root.
\end{itemize}
On Dowling trees the operation of inner orbit contraction and extension are defined analogously as in $\TG$ with the only difference that for every edge $t$ that connects two zero vertices we have $\orb (t)=\{t\}$. Therefore the Dowling trees form an abstract simplicial complex that we will denote by $\TNG$.
  }\\

\begin{figure}[htbp]\label{trees}
\begin{center}
\scalebox{1}{
\psfrag{a}{(a)}
\psfrag{b}{(b)}
\psfrag{x}{$\orb (\{0,1,2, \w , \ww\})$}
\psfrag{y}{$\orb (\{1,\ww\})$}
\psfrag{1}{$1$}
\psfrag{2}{$2$}
\psfrag{3}{$3$}
\psfrag{u}{$\w$}
\psfrag{d}{$\ww$}
\psfrag{t}{$\www$}
\psfrag{T}{$T_1$}
\psfrag{Q}{$T_2$}
\includegraphics[scale=0.6]{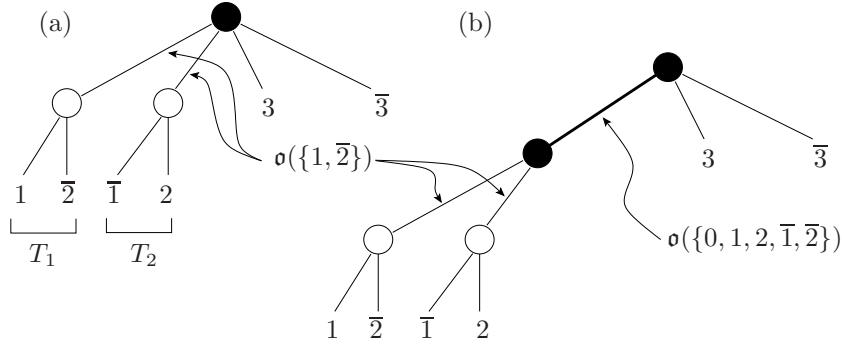}
} 
\caption{(a): the tree corresponding to $\{1,\ww\}$. (b): the Dowling tree constructed from the nested set $\{1\ww , 012\}$. The zero vertices are black.
{\em As above, only a representative of every nonsingleton orbit is indicated. See figure 1.}}
\end{center}
\end{figure}

We state the theorem relating Dowling trees and Nested set complexes
of Dowling lattices. The way of encoding trees with nested sets is the
same as in \cite{F}.

{\thm $\TNG$ is isomorphic to $\RN(\JJ^G, \Dl)$ as an abstract simplicial complex.}

{\pf Consider a simplex $X\in \RN(\JJ^G,\Dl)$. Since $\dl$ is an order
  ideal in $\Dl$, we may choose a linear extension of the ordering in
  $X$ such that all elements of type $1$ come before all those of type
  $0$. We will perform our inner orbit extensions according to the
  chosen linear order of $X$. After having exhausted all elements of
  type $1$ we are clearly left with a tree $T\in \TG$, that can be
  turned into a good Dowling tree just by declaring the root as the
  only zero vertex. On this tree we now have to perform the `type $0$' - orbit extensions. 

So let $\omega$ be a type-$0$ partition, with zero block $S_0$. As above, there is a vertex $v$ of $T$ such that the union of the leaves in the connected components $T_1,\dots , T_s$ of $T$ not containing the root that arise by deleting $v$ is minimal with the property of containing the set $S_0\setminus \{0\}$.
In particular, this $v$ is fixed by the action of $G$ and therefore is
a zero vertex; by \ref{minbuild}(4) we know that the type $0$ elements
of $X$
that corresponding to already performed extensions lie on a chain
below $\omega$, so that their zero blocks are all contained in $S_0$. Thus, $v$ can be only the root.

Let $\tau_1, \dots , \tau_k$ denote the elements of $X$ that are maximal among those below $\omega$. By construction, to every $\tau_i$ corresponds an inner orbit of edges that are incident to the root. Again by construction, all elements in $S_0\setminus\{0\}$ that are not contained in a nonsingleton block of some $\tau_i$ are directly appended to the root. 
We may then renumber the $T_i$'s in such a way that the union of the labels of the leafs of the first $s'$ trees is exactly $S_0\setminus \{0\}$. Note that $s'<s$ because $\hat{1}\not\in X$ (in figure 2 we have $s'=2$ and the corresponding trees $T_1,T_2$ are indicated).

 Then we build a tree starting with an edge $t$ that joins the root
 and a new vertex $w$ (which we declare to be a zero vertex). Below
 $w$ we grow $s'$ edges $e_1,\dots e_{s'}$ , and append to those the
 trees $T_1,\dots T_{s'}$. The trees $T_{s'+1},\dots T_s$ will be
 appended directly to the root via edges $e_{s'+1},\dots e_s$.

Now check that this is again a Dowling tree: we only have to worry about the zero vertices. Both the sets $\bigcup_{i\leq s'}\lambda_{e_i}$ and $\bigcup_{i>s'}\lambda_{e_i}$ contain the full orbit of each of their elements, and therefore property (1), (2) and (3) follow immediately. For property (4) recall Lemma \ref{minbuild} to see that all zero vertices (except the root) are in some $T_i$ with $i\leq s'$.

We have thus constructed a unique Dowling tree $T(X)$ from a nested
set $X\in \NN(\JJ^G,\Dl)$. The inverse operation is now easy: given a
Dowling tree $T$ identify the orbits of all inner edges under the
action of $G$, and note that the set of corresponding one-block-orbit
elements of $\Dl$ is nested.

\label{tng} It is clear that the bijection $T:\RN(\JJ^G,\Dl)\rightarrow \TNG$ extends to an isomorphism of simplicial complexes, if we take the operation of orbit contraction as boundary operator in $\TNG$.\hfill$\square$ }\\

Summarizing, we can formulate the following corollary, that is a suggestive counterpart of corollary \ref{hultman}.

{\crl\label{subD} The complex of Dowling trees $\TNG$ is a pure simplicial complex of dimension $(n-2)$. It is subdivided by the reduced order complex $\RD(\Dl)$ of the Dowling lattice. Their realizations are therefore PL-homeomorphic. They are homotopy equivalent to a wedge of $$(|G|+1)(2|G|+1)\dots((n-1)|G|+1)$$ spheres.}\\[1mm]
{\noindent{\bf Example.} For the examples considered above, where
  $n=3$ and $G=\mathbb{Z}_2$, we have that the complexes
  $\RD(\mathcal{Q}_3^0(\mathbb{Z}_2))$ and
  $\mathcal{T}_3^{\mathbb{Z}_2}$ are each homotopy equivalent of a wedge of $(2-1)(4-1)=3$ circles, while $\RD(\mathcal{Q}_3(\mathbb{Z}_2))$ and $\mathcal{T}_3(\mathbb{Z}_2)$ have the homotopy type of a wedge of $(2+1)(2\cdot 2 +1)=15$ circles.}\\

\subsection{From $\TG$ to $\TNG$} The description in terms of nested
  set complexes allows us to explicitly reconstruct $\TNG$ from $\TG$ by successively coning over subcomplexes having the homotopy type of wedges of $(n-3)$-spheres. This gives another proof of the fact that $\TG$ is Cohen-Macaulay and allows to explicitly calculate the difference of the number of spheres in the homotopy type of $\TG$ and $\TNG$.

First of all we want to distinguish three types of simplices in
$\TNG$. We call simplices of {\em type 0}, respectively of {\em type
  1}, those simplices consisting only of elements of type $0$,
respectively of type $1$. The nested sets containing elements of both
types will be called simplices {\em of mixed type}. We remark that the subcomplex given by the simplices of type $1$ is exactly $\TG$, and that any simplex $X$ of mixed type is contained in the star of a unique maximal simplex $X_0$ of type $0$, namely $X_0 = X\setminus \II^G$.

The idea is therefore to start with $\TG$ and glue successively the stars of all simplices of type $0$. Topologically, this means coning over the link of those simplices: to keep track of the change of topology, we need some definitions and a lemma.\\

{\df\label{def_k} Let $\TT_J$ denote the subcomplex of $\TNG$ consisting of all simplices of type $0$, i.e.\begin{center}
$\TT_J:=\{ X \in \TNG \vert X\subset \JJ^G \setminus \II^G\} $.\end{center}
From the above considerations we know that any $X\in \TT_J$ is a chain $\omega = \omega_1 < \omega_2 < \dots <\omega_\ell $ of elements of type $0$. The { length} of the chain is the number of its elements and will be denoted $\ell(\omega)=\ell(X)$. The associated partitions $\uo_i \in \Pno$ have only one nonsingleton block, namely the one containing $0$, which we call $w_i$. Setting $w_0 := \{0\}$ and $w_{\ell +1}:= \{0,1,2 ,\dots ,n\}$, we define numbers $p_0(\omega),\dots p_{\ell}(\omega) \in \mathbb{N}$ as\begin{center} $p_i(\omega):= \vert w_{i+1} \setminus w_{i}\vert $.  \end{center} 
If the chain is understood, we will simply write $p_i$.
For $m=1,\dots n-1$ we define the subcomplex of $\TNG$ consisting of $\TG$ and the stars of all simplices $X\in \TT_J$ with $\ell(X)\leq m$:
\begin{center}$\KK_m:= \TG \cup \{X\in \TNG \vert\, \vert X\cap \JJ^G\vert\leq m\}. $ \end{center}}

{\lm The link of any $X\in \TT_J$ with $\ell (X) = m $ in $\KK_m$ is 
$$ lk_{\KK_m}(X) \simeq \RD(B_m) \ast \RD(\mathcal{Q}^G_{p_0}) \ast \dots \ast \RD(\mathcal{Q}^G_{p_{m}}),  $$ where $B_m$ denotes the boolean lattice on $m$ elements.}

{\pf Any simplex $Y$ in the link can be written as $$Y= Y' \amalg Y_0 \amalg \dots Y_m,$$
where $Y'$ is a (proper!) subset of $X$, and $Y_i$ is a nested subset of $\II^G$ such that the only nonsingleton block of the associated partitions in $\Pno$ contains only elements from $w_{i+1}\setminus w_i$. The subcomplex of such $Y_i$ can of course be identified with $\RN(\II^G, \mathcal{Q}_{p_i}^G)$, whereas the possible choices of $Y'$ give a subcomplex with a face lattice that can be identified with the proper part of $B_m$, the boolean lattice on $m$ elements. Note that any choice of $Y'\in B_m\setminus \{\hat{1}\}$ and $Y_i \in\RN(\II^G, \mathcal{Q}_{p_i}^G)$ gives a simplex in the link. 

Since all complexes $\RD(B_m)$ and $\RN(\II^G, \mathcal{Q}^G_{p_i})$ are flag complexes, we have that the link of $X$ is a simplicial complex that is isomorphic to the join
$$ \RD(B_m) \ast \RN(\II^G, \mathcal{Q}_{p_0}^G)\ast\dots \ast \RN(\II^G, \mathcal{Q}_{p_m}^G).$$ With corollary \ref{subd_crl} the claim follows. \hfill$\square$ }\\

In order to simplify notation, let us define numbers $q_i^\omega$ associated to any chain $\omega$ that gives rise to a simplex in $\TT_J$. Recall definition \ref{def_k} and let $$q_i^\omega:= \prod_{j=1}^{p_i(\omega)-1} (j|G|-1).$$
The numbers $Q(\omega)$ are then defined for any chain $\omega$ as $$Q(\omega):= q_0^\omega q_1^\omega\dots q_{\ell(\omega)}^\omega.$$

With these definitions we can state the theorem, which now follows
easily from our previous work.

{\thm\label{nostro} The link of any $X\in \TT_J$ with $\ell (X) = m $ in $\KK_m$ is homotopy equivalent to a wedge of $Q(\omega)$ spheres of dimension $(n-3)$, where $\omega$ is the chain obtained by ordering the elements of $X$. Each of those spheres is the boundary of an $(n-2)$-ball in $\KK_m$. }

{\pf After Hultman \cite{Hu} we know that, for any $p_i(\omega)$, $\RD(\mathcal{Q}_{p_i(\omega)}^G)$ is homotopy equivalent to a wedge of $q_i^\omega$ spheres of dimension $(p_i(\omega)-2)$. It is a standard fact that $\RD(B_m)\simeq S^{(m-1)}$. We have then to compute the homotopy type of
$$S^{(m-1)}\ast \bigvee_{q_0^\omega}S^{(p_0 - 2)} \ast \dots \ast\bigvee_{q_m^\omega}S^{(p_m - 2)},  $$
where the index under the wedges indicates how many copies of the
corresponding sphere come into the game. By basic topological facts we
may rewrite this as:
$$ S^{(m-2)} \ast \bigvee_{q_0^\omega \dots q_m^\omega} S^{\sum_{i=0}^m (p_i - 2) + m} = \bigvee_{Q(\omega)} S^{m-2 + n -2(m+1) +m +1} = \bigvee_{Q(\omega)} S^{(m-3)}, $$
where in the second equality we used that $p_0 + p_1 + \dots p_m = n$. This proves the first part of the corollary.

The last assertion is proved by induction on $m$, after remarking that actually the link of $X$ in $\KK_m$ is contained in $\KK_{m-1}$ (we define $\KK_0=\TG$). For $m=1$ the assertion holds because $\TG$ is CM of dimension $(n-2)$, thus each $(n-3)$-sphere bounds a ball. Let the claim hold for $m\geq 1$. Then in particular $\KK_m$ was obtained by repeatedly coning over spheres that were already boundaries - therefore $\KK_m$ is also CM of dimension $(n-2)$, and any of its $(n-3)$-cycles bounds.\hfill$\square$}\\

{\noindent\bf Note:} Since $p_i(\omega) < n$, we need the result of \cite{Hu} only in dimension strictly smaller than the one in which the conclusion of the corollary holds. Therefore we may in principle omit the use of \cite{Hu}, thus reproving fully independently the result, by an induction on $n$.

{\rem\label{ultimo} We have proved that any chain $\omega\in \RD(B_n)$ indexes a simplex of $\TT_J$ that contributes $Q(\omega)$ times to the difference of the number of spheres between the homotopy types of $\TG$ and $\TNG$.}\\

{\noindent\bf Example:} For our favourite example $\mathcal{Q}_3(\mathbb{Z}_2)$, we have $12$ chains in $\RD(B_3)$, each with $Q(\omega)=1$, therefore $\sum_{\omega\in \RD(B_3)}Q(\omega)=12$, which in fact gives $12+3 = 15$.\\

We may even combine the results of Dowling about $\RD(\Dl)$, of
Hultman about $\TG$ and our above considerations to state the following arithmetic equality:

{\crl\label{numerology} Let integers $k\geq 1$ and $n\geq 2$ be given, and  for $\pi\in\Pi_n$ let $h(\pi, j)$ denote the height of the $j$-th column of the Young tableau of $\pi$. Then 
$$ \prod_{j=1}^{n}(jk+1) - \prod_{j=1}^{n}(jk-1)= \sum_{\sigma \in \Pi_n} \prod_{j=1}^{n}(jk-1)^{h(\sigma, j)}.$$}


\begin{thebibliography}{BLSZ}


\bibitem{AK} F. Ardila, C. Klivans; {\em The Bergman complex of a
    matroid and phylogenetic trees.}  J. Combin. Theory Ser. B  {\bf 96}  (2006),  no. 1, 38-49.


\bibitem{BHV}L.\ J.\ Billera, S. \ P. \ Holmes, K. \ Vogtmann; {\em Geometry of the space of phylogenetic trees.}  Adv. in Appl. Math.  {\bf 27 } (2001),  no. 4, 733-767.

\bibitem{Bj} A. Bj\"orner; {\em Shellable and Cohen-Macaulay partially ordered sets.} Trans. Amer. Math. Soc. {\bf 260} (1980), 159-183.







\bibitem{Bo} J. M. Boardman; {\em Homotopy structures and the language of trees.} Proc. Symp. Pure Math. {\bf 22} (1971), 37-58.



\bibitem{CD} S. \v Cuki\'c, E. Delucchi; {\em Shellable simplicial spheres via combinatorial blowups.} ArXiv \texttt{math.CO/0602101}.




\bibitem{D} E. Delucchi; {\em Subdivision of complexes of $k$-trees.} ArXiv \texttt{math.CO/0509378}.


\bibitem{Dow1} T. A. Dowling; {\em A $q$-analog of the partition lattice.} A survey of combinatorial theory (Proc. Internat. Sympos., Colorado State Univ., 1971), pp. 101-115. 

\bibitem{Dow}  T. A. Dowling; {\em A class of geometric lattices based on finite groups.}  J. Combinatorial Theory Ser. B  {\bf 14}  (1973), 61-86. (Erratum:  J. Combinatorial Theory Ser. B  {\bf 15}  (1973), 211.)

\bibitem{ER} R. Ehrenborg, M. A. Readdy; {\em The Dowling transform of subspace arrangements.}
J. Combin. Theory Ser. A {\bf 91} (2000), no. 1-2, 322-333.

\bibitem{ER2} R. Ehrenborg, M. A. Readdy; {\em
On flag vectors, the Dowling lattice, and braid arrangements.}  Discrete Comput. Geom.  {\bf 21}  (1999),  no. 3, 389-403.


\bibitem{F} E. M. Feichtner; {\em Complexes of trees and nested set
    complexes.} ArXiv math.CO/0409235, to appear in Pacific J. of Math. 

\bibitem{FK} E. M. Feichtner, D. N. Kozlov; {\em Incidence combinatorics of resolutions.}  Selecta Math. (N.S.)  {\bf 10}  (2004),  no. 1, 37-60.

\bibitem{FK2} E. M. Feichtner, D.  N. Kozlov; {\em Abelianizing the real permutation action via blowups.}  Int. Math. Res. Not.  (2003),  no. 32, 1755-1784.

\bibitem{FM} E. M. Feichtner, I. M\"uller; {\em On the topology of nested set complexes.} Proc. Amer. Math. Soc. {\bf 133} (2005), no.4, 999-1006.

\bibitem{FS} E. M.  Feichtner, B.  Sturmfels; {\em Matroid polytopes, nested sets, and Bergman fans.} ArXiv \texttt{math.CO/0411260}, to appear in Port. Math. (N.S.).


\bibitem{FY} E. M. Feichtner, S. Yuzvinsky; {\em Chow rings of toric varieties defined by atomic lattices.}  Invent. Math.  {\bf 155}  (2004),  no. 3, 515-536.

\bibitem{Han1} P. Hanlon; {\em Otter's method and the homology of homeomorphically irreducible $k$-trees.}  J. Combin. Theory Ser. A  {\bf 74}  (1996),  no. 2, 301-320.


\bibitem{Han} P. Hanlon; {\em The generalized Dowling lattices.}  Trans. Amer. Math. Soc.  {\bf 325}  (1991),  no. 1, 1-37.  

\bibitem{HW} P. Hanlon, M. Wachs; {\em On Lie $k$-algebras.}  Adv. Math.  {\bf 113 } (1995),  no. 2, 206-236.

\bibitem{Hu} A. Hultman; {\em The topology of spaces of phylogenetic trees with symmetry.} Preprint available at \texttt{http://www.math.kth.se/\~{}hultman/}.


\bibitem{RW} A. Robinson, S. Whitehouse; {\em The tree representation of $\Sigma\sb {n+1}$. } J. Pure Appl. Algebra  {\bf 111}  (1996),  no. 1-3, 245-253.


\bibitem{Sp} E. H. Spanier; {\em Algebraic topology.} Springer-Verlag, New York-Berlin, 1981. 

\bibitem{Sta} R. P.  Stanley; {\em Supersolvable lattices.}  Algebra Universalis  {\bf 2}  (1972), 197--217. 

\bibitem{TZ} H. Trappmann, G. M. Ziegler; {\em Shellability of complexes of trees.} J. Combin. Theory Ser. A {\bf 82} (1998), 168-178.

\bibitem{V} K. Vogtmann; {\em Local structure of some ${\rm Out}(F\sb n)$-complexes.}  Proc. Edinburgh Math. Soc. (2)  {\bf 33}  (1990),  no. 3, 367-379.


\end{thebibliography}
\end{document}